# New proofs of the duplication and multiplication formulae for the gamma and the Barnes double gamma functions

Donal F. Connon

dconnon@btopenworld.com

26 March 2009

**Abstract**

New proofs of the duplication formulae for the gamma and the Barnes double gamma functions are derived using the Hurwitz zeta function. Concise derivations of Gauss's multiplication theorem for the gamma function and a corresponding one for the double gamma function are also reported. This paper also refers to some connections with the Stieltjes constants.

## 1. Legendre's duplication formula for the gamma function

Hansen and Patrick [11] showed in 1962 that the Hurwitz zeta function could be written as

(1.1) $$\varsigma(s,x) = 2^s \varsigma(s,2x) - \varsigma\left(s, x+\frac{1}{2}\right)$$

and, by analytic continuation, this holds for all $s$. Differentiation results in

(1.2) $$\varsigma'(s,x) = 2^s \varsigma'(s,2x) + 2^s \log 2\, \varsigma(s,2x) - \varsigma'\left(s, x+\frac{1}{2}\right)$$

and with $s = 0$ we have

(1.3) $$\varsigma'(0,x) = \varsigma'(0,2x) + \log 2\, \varsigma(0,2x) - \varsigma'\left(0, x+\frac{1}{2}\right)$$

We recall Lerch's identity for $\mathrm{Re}(s) > 0$

(1.4) $$\log \Gamma(x) = \varsigma'(0,x) - \varsigma'(0) = \varsigma'(0,x) + \frac{1}{2}\log(2\pi)$$

The above relationship between the gamma function and the Hurwitz zeta function was established by Lerch in 1894 (see, for example, Berndt's paper [6]). A different proof is contained in [9].

We have the well known relationship between the Hurwitz zeta function and the Bernoulli polynomials $B_n(u)$ (for example, see Apostol's book [4, pp. 264-266]).

(1.5) $$\varsigma(-m, x) = -\frac{B_{m+1}(x)}{m+1} \quad \text{for } m \in \mathbf{N}_o$$

which gives us the well-known formula

$$\varsigma(0, x) = \frac{1}{2} - x$$

Therefore we have from (1.3) and (1.4)

(1.6) $$\log \Gamma(x) + \log \Gamma\left(x + \frac{1}{2}\right) = \log \Gamma(2x) + (1 - 2x) \log 2 + \log \sqrt{\pi}$$

and hence we obtain Legendre's duplication formula [16, p.240] for the gamma function

(1.7) $$2^{2x-1} \Gamma(x) \Gamma\left(x + \frac{1}{2}\right) = \sqrt{\pi}\, \Gamma(2x)$$

Hansen and Patrick [11] also showed that

(1.8) $$\sum_{r=1}^{q-1} \varsigma\left(s, \frac{r}{q}\right) = (q^s - 1)\varsigma(s)$$

Differentiation results in

(1.9) $$\sum_{r=1}^{q-1} \varsigma'\left(s, \frac{r}{q}\right) = (q^s - 1)\varsigma'(s) + q^s \varsigma(s) \log q$$

and with $s = 0$ we have

(1.10) $$\sum_{r=1}^{q-1} \varsigma'\left(0, \frac{r}{q}\right) = -\frac{1}{2} \log q$$

Substituting Lerch's identity (1.4) we get

(1.11) $$\sum_{r=1}^{q-1} \log \Gamma\left(\frac{r}{q}\right) = \frac{q-1}{2} \log(2\pi) - \frac{1}{2} \log q$$

and with $q = 2$ this immediately gives us the well-known result [15, p.3]

(1.12) $$\log \Gamma\left(\frac{1}{2}\right) = \frac{1}{2}\log \pi$$

It should be noted that the proof of the above identity is dependent on Lerch's identity which may be derived without assuming any prior knowledge of (1.12). In the author's view, this is akin to the marvel experienced when first confronted with a derivation of Euler's integral

$$\int_0^{\pi/2} \log \sin x \, dx = -\frac{\pi}{2}\log 2$$

With $q = 4$ in (1.10) we see that

$$\log \Gamma\left(\frac{1}{4}\right) + \log \Gamma\left(\frac{1}{2}\right) + \log \Gamma\left(\frac{3}{4}\right) = \frac{3}{2}\log(2\pi) - \frac{1}{2}\log 4$$

and thus

(1.13) $$\log \Gamma\left(\frac{1}{4}\right) + \log \Gamma\left(\frac{3}{4}\right) = \log \pi + \frac{1}{2}\log 2$$

which of course may also be easily obtained directly from Euler's reflection formula for the gamma function.

With $s = -1$ in (1.9), and using $\varsigma(-1) = -\frac{1}{12}$, we obtain

(1.14) $$\sum_{r=1}^{q-1} \varsigma'\left(-1, \frac{r}{q}\right) = (q^{-1} - 1)\varsigma'(-1) - \frac{1}{12q}\log q$$

and with $q = 2$ we have

(1.15) $$\varsigma'\left(-1, \frac{1}{2}\right) = -\frac{1}{2}\varsigma'(-1) - \frac{1}{24}\log 2$$

which we shall also see below in (3.2).

## 2. Gauss's multiplication theorem for the gamma function

The general Kubert identity is derived in [13, p.169]

(2.1) $$\Phi(s, x, z) = q^{-s} \sum_{r=0}^{q-1} \Phi\left(s, \frac{r+x}{q}, z^q\right)$$

where $\Phi(s, x, z)$ is the Hurwitz-Lerch zeta function

$$(2.2) \qquad \Phi(s, x, z) = \sum_{n=0}^{\infty} \frac{z^n}{(n+x)^n}$$

We see that $\Phi(s, x, 1) = \varsigma(s, x)$ and therefore we have

$$(2.3) \qquad q^s \varsigma(s, x) = \sum_{r=0}^{q-1} \varsigma\left(s, \frac{r+x}{q}\right)$$

which corresponds with (1.8) when $x = 1$.

Differentiation results in

$$(2.4) \qquad q^s \varsigma'(s, x) + q^s \varsigma(s, x) \log q = \sum_{r=0}^{q-1} \varsigma'\left(s, \frac{r+x}{q}\right)$$

and letting $s = 0$ and substituting Lerch's identity (1.4) we get

$$(2.5) \qquad \log \Gamma(x) = \sum_{r=0}^{q-1} \log \Gamma\left(\frac{r+x}{q}\right) - \frac{(q-1)}{2} \log(2\pi) - \left(\frac{1}{2} - x\right) \log q$$

or

$$(2.6) \qquad (2\pi)^{(q-1)/2} n^{(1/2)-x} \Gamma(x) = \prod_{r=0}^{q-1} \Gamma((r+x)/q)$$

which is Gauss's multiplication theorem for the gamma function [3, p.23]. I subsequently discovered that a similar procedure was employed in Milnor's paper [14].

Letting $s = 1 - n$ in (2.3) gives us

$$\varsigma(1-n, x) = \frac{1}{q^{n-1}} \sum_{r=0}^{q-1} \varsigma\left(1-n, \frac{r+x}{q}\right)$$

and using (1.5) results in

$$B_n(x) = \frac{1}{q^{n-1}} \sum_{r=0}^{q-1} B_n\left(\frac{r+x}{q}\right)$$

where the substitution $x \to qx$ gives us the multiplication formula for the Bernoulli polynomials [15, p.60]

$$B_n(qx) = \frac{1}{q^{n-1}} \sum_{r=0}^{q-1} B_n\left(x + \frac{r}{q}\right)$$

Differentiation of (2.5) gives us [15, p.12]

(2.7) $$\psi(x) = \log q + \frac{1}{q} \sum_{r=0}^{q-1} \psi\left(\frac{r+x}{q}\right)$$

and further differentiations give us

$$\psi^{(n)}(x) = \frac{1}{q^{n+1}} \sum_{r=0}^{q-1} \psi^{(n)}\left(\frac{r+x}{q}\right)$$

Since [15, p.22]

$$\psi^{(n)}(x) = (-1)^{n+1} n! \varsigma(n+1, x)$$

we see that this results in

$$\varsigma(n+1, x) = \frac{1}{q^{n+1}} \sum_{r=0}^{q-1} \varsigma\left(n+1, \frac{r+x}{q}\right)$$

which is a particular case of (2.3) for positive integer values of $s$.

Hansen and Patrick [11] also showed that

(2.9) $$\sum_{r=1}^{q} \varsigma\left(s, \frac{r}{q} - b\right) = q^s \varsigma(s, 1 - bq)$$

and letting $b = -\frac{x}{q}$ we have

$$\sum_{r=1}^{q} \varsigma\left(s, \frac{r+x}{q}\right) = q^s \varsigma(s, 1+x)$$

Noting that

(2.10) $$\varsigma(s, 1+x) = \varsigma(s, x) - x^{-s}$$

this becomes

$$\sum_{r=1}^{q} \varsigma\left(s, \frac{r+x}{q}\right) = q^s \varsigma(s,x) - \frac{q^s}{x^s}$$

which may be written as

$$\sum_{r=0}^{q-1} \varsigma\left(s, \frac{r+x}{q}\right) - \varsigma\left(s, \frac{x}{q}\right) + \varsigma\left(s, 1+\frac{x}{q}\right) = q^s \varsigma(s,x) - \frac{q^s}{x^s}$$

Letting $x \to \frac{x}{q}$ in (2.10) we then obtain another derivation of (2.3).

## 3. Duplication formula for the Barnes double gamma function

With $s = -1$ in (1.2) we have

(3.1) $\quad \varsigma'(-1,x) = \frac{1}{2}\varsigma'(-1,2x) + \frac{1}{2}\log 2\, \varsigma(-1,2x) - \varsigma'\left(-1, x+\frac{1}{2}\right)$

and using (1.5) we have

$$\varsigma(-1,2x) = -\frac{1}{2}\left(4x^2 - 2x + \frac{1}{6}\right)$$

For example, equation (3.1) also gives us for $x = 1/2$

(3.2) $\quad \varsigma'\left(-1, \frac{1}{2}\right) = -\frac{1}{24}\log 2 - \frac{1}{2}\varsigma'(-1)$

We have the Gosper/Vardi functional equation for the Barnes double gamma

(3.3) $\quad \varsigma'(-1,x) = \varsigma'(-1) - \log G(1+x) + x\log \Gamma(x)$

which was derived by Vardi in 1988 and also by Gosper in 1997 (see [1]). A different derivation is given in equation (4.3.126) of [9].

Using this and (3.3) we may easily deduce that

(3.4) $\quad \log G\left(\frac{1}{2}\right) = -\frac{1}{4}\log \pi + \frac{1}{24}\log 2 + \frac{3}{2}\varsigma'(-1)$

as originally determined by Barnes [5] in 1899.

Combining (3.1) and (3.3) results in

$$-\log G(1+x) + x\log \Gamma(x) = -\frac{1}{2}\log G(1+2x) + x\log \Gamma(2x)$$

$$-\frac{1}{4}\left(4x^2 - 2x + \frac{1}{6}\right)\log 2 - \frac{3}{2}\varsigma'(-1) - \log G\left(x + \frac{3}{2}\right) + \left(x + \frac{1}{2}\right)\log \Gamma\left(x + \frac{1}{2}\right)$$

Since $G(1+x) = G(x)\Gamma(x)$ this may be written as

$$-\log G(x) - \log \Gamma(x) + x\log \Gamma(x) = -\frac{1}{2}\log G(2x) - \log \Gamma(2x) + x\log \Gamma(2x)$$

$$-\frac{1}{4}\left(4x^2 - 2x + \frac{1}{6}\right)\log 2 - \frac{3}{2}\varsigma'(-1) - \log G\left(x + \frac{1}{2}\right) - \log \Gamma\left(x + \frac{1}{2}\right) + \left(x + \frac{1}{2}\right)\log \Gamma\left(x + \frac{1}{2}\right)$$

and using (1.6) we thereby obtain the duplication formula for the Barnes double gamma function. In 1899 Barnes developed a multiplication formula for $G(nx)$ (see [15, p.30]) and a particular case is set out below [15, p.29]

(3.5) $$G^2(x)G^2\left(x + \frac{1}{2}\right)\Gamma(x) = J(x)G(2x)$$

where for convenience $J(x)$ is defined by

$$\log J(x) = \frac{1}{4} - 3\log A + \left(-2x^2 + 3x - \frac{11}{12}\right)\log 2 + \left(x - \frac{1}{2}\right)\log \pi$$

A different derivation of this duplication formula was given by Choi [7] in 1996 where he used the double Hurwitz zeta function defined by

$$\varsigma_2(s,a) = \sum_{k_1,k_2 \geq 0} (a + k_1 + k_2)^{-s}$$

**4. A multiplication formula for the Barnes double gamma function**

With $s = -1$ in (2.4) we have

$$\varsigma'(-1,x) - \frac{1}{2}B_2(x)\log q = q\sum_{r=0}^{q-1}\varsigma'\left(-1, \frac{r+x}{q}\right)$$

and with the Gosper/Vardi functional equation (3.3) this becomes

$$\varsigma'(-1) - \log G(1+x) + x \log \Gamma(x) - \frac{1}{2} B_2(x) \log q$$

$$= q^2 \varsigma'(-1) - q \sum_{r=0}^{q-1} \log G\left(1 + \frac{r+x}{q}\right) + q \sum_{r=0}^{q-1} \left(\frac{r+x}{q}\right) \log \Gamma\left(\frac{r+x}{q}\right)$$

However, it is not immediately clear how this may be expressed in the form of the multiplication formula originally derived by Barnes [5, p.291].

Substituting $x = qt$ we have

$$\varsigma'(-1) - \log G(1+qt) + qt \log \Gamma(qt) - \frac{1}{2} B_2(qt) \log q$$

$$= q^2 \varsigma'(-1) - q \sum_{r=0}^{q-1} \log G\left(1 + t + \frac{r}{q}\right) + q \sum_{r=0}^{q-1} \left(t + \frac{r}{q}\right) \log \Gamma\left(t + \frac{r}{q}\right)$$

**5. Other multiple gamma functions**

Adamchik [2] has shown that for $\mathrm{Re}(x) > 0$

(5.1) $\qquad \varsigma'(-n, x) - \varsigma'(-n) = (-1)^n \sum_{k=0}^{n} k! Q_{k,n}(x) \log \Gamma_{k+1}(x)$

where the polynomials $Q_{k,n}(x)$ are defined by

$$Q_{k,n}(x) = \sum_{j=k}^{n} (1-x)^{n-j} \binom{n}{j} \begin{Bmatrix} j \\ k \end{Bmatrix}$$

and $\begin{Bmatrix} j \\ k \end{Bmatrix}$ are the Stirling subset numbers defined by

$$\begin{Bmatrix} j \\ k \end{Bmatrix} = k \begin{Bmatrix} n-1 \\ k \end{Bmatrix} + \begin{Bmatrix} n-1 \\ k-1 \end{Bmatrix}, \quad \begin{Bmatrix} n \\ 0 \end{Bmatrix} = \begin{cases} 1, & n=0 \\ 0, & n \neq 0 \end{cases}$$

We have [15, p.39]

(5.2) $\qquad G_n(x+1) = G_n(x) G_{n-1}(x) \qquad \Gamma_n(x) = [G_n(x)]^{(-1)^{n-1}}$

and it is easily seen that

$$\log G_n(x+1) = \log G_n(x) + \log G_{n-1}(x)$$

$$\log \Gamma_n(x) = (-1)^{n-1} \log G_n(x)$$

and from this we obtain

(5.3) $$\Gamma_n(x+1) = \frac{\Gamma_n(x)}{\Gamma_{n-1}(x)}$$

Particular cases of (5.1) are

(5.4) $$\varsigma'(-1,x) - \varsigma'(-1) = x \log \Gamma(x) - \log G(x+1)$$

(5.5) $$\varsigma'(-2,x) - \varsigma'(-2) = 2\log \Gamma_3(x) + (3-2x)\log G(x) - (1-x)^2 \log \Gamma(x)$$

Hence, using (2.3) we may obtain multiplication formulae for the higher order multiple gamma functions.

## 6. Some connections with the Stieltjes constants

The generalised Euler-Mascheroni constants $\gamma_n$ (or Stieltjes constants) are the coefficients of the Laurent expansion of the Riemann zeta function $\varsigma(s)$ about $s=1$

(6.1) $$\varsigma(s) = \frac{1}{s-1} + \sum_{n=0}^{\infty} \frac{(-1)^n}{n!} \gamma_n (s-1)^n$$

The Stieltjes constants $\gamma_n(x)$ are the coefficients in the Laurent expansion of the Hurwitz zeta function $\varsigma(s,u)$ about $s=1$

(6.2) $$\varsigma(s,x) = \sum_{n=0}^{\infty} \frac{1}{(n+x)^s} = \frac{1}{s-1} + \sum_{n=0}^{\infty} \frac{(-1)^n}{n!} \gamma_n(x)(s-1)^n$$

and $\gamma_0(x) = -\psi(x)$, where $\psi(x)$ is the digamma function which is the logarithmic derivative of the gamma function $\psi(x) = \frac{d}{dx} \log \Gamma(x)$. It is easily seen from the definition of the Hurwitz zeta function that $\varsigma(s,1) = \varsigma(s)$ and accordingly that $\gamma_n(1) = \gamma_n$.

Since $\lim_{s \to 1}\left[\varsigma(s) - \frac{1}{s-1}\right] = \gamma$ it is clear that $\gamma_0 = \gamma$. It may be shown, as in [12, p.4], that

(6.3) $$\gamma_n = \lim_{N \to \infty} \left[ \sum_{k=1}^{N} \frac{\log^n k}{k} - \frac{\log^{n+1} N}{n+1} \right] = \lim_{N \to \infty} \left[ \sum_{k=1}^{N} \frac{\log^n k}{k} - \int_1^N \frac{\log^n t}{t} dt \right]$$

where, throughout this paper, we define $\log^0 1 = 1$.

It was previously shown in [10] that

(6.4) $$\gamma_n(x) = -\frac{1}{n+1} \sum_{i=0}^{\infty} \frac{1}{i+1} \sum_{j=0}^{i} \binom{i}{j} (-1)^j \log^{n+1}(x+j)$$

We see from (6.2) that for $n \geq 0$

(6.5) $$\left. \frac{d^{n+1}}{ds^{n+1}} [(s-1)\varsigma(s,x)] \right|_{s=1} = (-1)^n (n+1) \gamma_n(x)$$

We multiply (2.3) by $(s-1)$

$$q^s (s-1) \varsigma(s,x) = \sum_{r=0}^{q-1} (s-1) \varsigma\left(s, \frac{r+x}{q}\right)$$

and, using the Leibniz rule to differentiate this $n+1$ times, we obtain

$$\sum_{k=0}^{n+1} \binom{n+1}{k} \frac{d^{n+1-k}}{ds^{n+1-k}} [(s-1)\varsigma(s,x)] \frac{d^k}{ds^k} q^s = \sum_{r=0}^{q-1} \frac{d^{n+1}}{ds^{n+1}} \left[ (s-1) \varsigma\left(s, \frac{r+x}{q}\right) \right]$$

Evaluating this at $s=1$ results in

$$q \log^{n+1} q + q \sum_{k=0}^{n} \binom{n+1}{k} (-1)^{n-k} (n-k+1) \gamma_{n-k}\left(\frac{r+x}{q}\right) \log^k q = \sum_{r=0}^{q-1} (-1)^n (n+1) \gamma_n\left(\frac{r+x}{q}\right)$$

(where we have isolated the $(n+1)$th term using $\lim_{s \to 1} [(s-1)\varsigma(s,x)] = 1$)

Using the binomial identity $\frac{n-k+1}{n+1} \binom{n+1}{k} = \binom{n}{k}$ this may be expressed as

(6.6) $$\sum_{r=0}^{q-1} \gamma_n\left(\frac{r+x}{q}\right) = q(-1)^n \frac{\log^{n+1} q}{n+1} + q \sum_{k=0}^{n} \binom{n}{k} (-1)^k \gamma_{n-k}(x) \log^k q$$

and noting that

$$\sum_{r=0}^{q-1} f\left(\frac{r+x}{q}\right) = \sum_{m=1}^{q} f\left(\frac{m-1+x}{q}\right) = \sum_{m=1}^{q-1} f\left(\frac{m-1+x}{q}\right) + f\left(\frac{q-1+x}{q}\right)$$

we see that for integers $q \geq 2$ and $x = 1$

(6.7) $$\sum_{r=1}^{q-1} \gamma_n\left(\frac{r}{q}\right) = -\gamma_n + q(-1)^n \frac{\log^{n+1} q}{n+1} + q \sum_{k=0}^{n} \binom{n}{k}(-1)^j \gamma_{n-k} \log^k q$$

which was previously derived by Coffey [8] using the relation (2.3).

With $n = 0$ in (6.7) we have

(6.8) $$\sum_{r=1}^{q-1} \gamma_0\left(\frac{r}{q}\right) = -\gamma + q\log q + q\gamma$$

Since $\psi(x) = -\gamma_0(x)$ we see from (2.7) that

(6.9) $$\gamma_0(x) = -\log q + \frac{1}{q}\sum_{r=0}^{q-1} \gamma_0\left(\frac{r+x}{q}\right)$$

and this concurs with (6.6) when $n = 0$.

With $x = 1$ this becomes

$$\gamma_0 = -\log q + \frac{1}{q}\sum_{r=0}^{q-1} \gamma_0\left(\frac{r+1}{q}\right) = -\log q + \frac{1}{q}\sum_{r=0}^{q-2} \gamma_0\left(\frac{r+1}{q}\right) + \frac{1}{q}\gamma_0$$

and therefore we obtain (6.8) again.

Letting $x \to 1+x$ in (6.9) we obtain

$$\gamma_0(1+x) = -\log q + \frac{1}{q}\sum_{r=0}^{q-1} \gamma_0\left(\frac{r+1+x}{q}\right)$$

$$= -\log q + \frac{1}{q}\sum_{m=1}^{q} \gamma_0\left(\frac{m+x}{q}\right)$$

and we then have

$$\gamma_0(1+x) = -\log q + \frac{1}{q}\sum_{m=0}^{q-1}\gamma_0\left(\frac{m+x}{q}\right) - \frac{1}{q}\gamma_0\left(\frac{x}{q}\right) + \frac{1}{q}\gamma_0\left(\frac{q+x}{q}\right)$$

Comparing this with (6.9) we obtain

(6.10) $\quad \gamma_0(1+x) = \gamma_0(x) - \frac{1}{q}\gamma_0\left(\frac{x}{q}\right) + \frac{1}{q}\gamma_0\left(1+\frac{x}{q}\right)$

For example, letting $q = 2$ we see that

$$\gamma_0(1+x) = \gamma_0(x) - \frac{1}{2}\gamma_0\left(\frac{x}{2}\right) + \frac{1}{2}\gamma_0\left(1+\frac{x}{2}\right)$$

Since $\psi(x) = -\gamma_0(x)$ we may express (6.10) as

$$\psi(1+x) - \psi(x) = \frac{1}{q}\psi\left(1+\frac{x}{q}\right) - \frac{1}{q}\psi\left(\frac{x}{q}\right)$$

and this may be easily verified by noting that [15, p.14]

$$\psi(1+x) - \psi(x) = \frac{1}{x} = \frac{1}{q}\frac{1}{(x/q)}$$

Letting $x \to 1+x$ in (6.6) we obtain

$$\sum_{r=0}^{q-1}\gamma_n\left(\frac{r+1+x}{q}\right) = q(-1)^n\frac{\log^{n+1}q}{n+1} + q\sum_{k=0}^{n}\binom{n}{k}(-1)^k\gamma_{n-k}(1+x)\log^k q$$

and noting that

$$\sum_{r=0}^{q-1}f\left(\frac{r+1+x}{q}\right) = \sum_{r=0}^{q-1}f\left(\frac{r+x}{q}\right) + f\left(1+\frac{x}{q}\right) - f\left(\frac{x}{q}\right)$$

we deduce that

$$q\sum_{k=0}^{n}\binom{n}{k}(-1)^k\gamma_{n-k}(x)\log^k q + \gamma_n\left(1+\frac{x}{q}\right) - \gamma_n\left(\frac{x}{q}\right) = q\sum_{k=0}^{n}\binom{n}{k}(-1)^k\gamma_{n-k}(1+x)\log^k q$$

or equivalently

$$\frac{1}{q}\left[\gamma_n\left(1+\frac{x}{q}\right)-\gamma_n\left(\frac{x}{q}\right)\right]=\sum_{k=0}^{n}\binom{n}{k}(-1)^k[\gamma_{n-k}(1+x)-\gamma_{n-k}(x)]\log^k q$$

With the reindexing $k=n-m$ we have

$$\sum_{k=0}^{n}\binom{n}{k}(-1)^k[\gamma_{n-k}(1+x)-\gamma_{n-k}(x)]\log^k q$$

$$=\sum_{m=n}^{0}\binom{n}{n-m}(-1)^{n-m}[\gamma_m(1+x)-\gamma_m(x)]\log^{n-m} q$$

$$=(-1)^n\log^n q\sum_{m=0}^{n}\binom{n}{m}(-1)^m[\gamma_m(1+x)-\gamma_m(x)]\log^{-m} q$$

arXiv:0710.4023 [pdf]

[10] D.F. Connon, Some series and integrals involving the Riemann zeta function, binomial coefficients and the harmonic numbers. Volume II(b), 2007.
arXiv:0710.4024 [pdf]

[11] E.R. Hansen and M.L. Patrick, Some Relations and Values for the Generalized Riemann Zeta Function.
Math. Comput., Vol. 16, No. 79. (1962), pp. 265-274.

[12] A. Ivić, The Riemann Zeta- Function: Theory and Applications.
Dover Publications Inc, 2003.

[13] S. Kanemitsu and H. Tsukada, Vistas of special functions.
World Scientific Publishing Co. Pte. Ltd., 2007.

[14] J. Milnor, On polylogarithms, Hurwitz zeta functions, and the Kubert identities.
L'Enseignement Math., 29 (1983), 281-322.

[15] H.M. Srivastava and J. Choi, Series Associated with the Zeta and Related Functions. Kluwer Academic Publishers, Dordrecht, the Netherlands, 2001.

[16] E.T. Whittaker and G.N. Watson, A Course of Modern Analysis: An Introduction to the General Theory of Infinite Processes and of Analytic Functions; With an Account of the Principal Transcendental Functions. Fourth Ed., Cambridge University Press, Cambridge, London and New York, 1963.

Donal F. Connon
Elmhurst
Dundle Road
Matfield
Kent TN12 7HD
dconnon@btopenworld.com
14